\documentclass[11pt,a4paper]{article}

\usepackage[UTF8]{ctex}
\usepackage[a4paper,margin=1in]{geometry}

\usepackage{amsmath}
\usepackage{amssymb}
\usepackage{amsthm}
\usepackage{mathtools}
\usepackage{mathrsfs}
\usepackage{bm}

\usepackage{graphicx}
\usepackage{booktabs}
\usepackage{enumitem}

\usepackage[colorlinks=true,
            linkcolor=blue,
            citecolor=blue,
            urlcolor=blue]{hyperref}

\providecommand{\AuthorMark}[1]{}
\providecommand{\TitleMark}[1]{}

\providecommand{\Keywords}[1]{%
  \par\noindent\textbf{Key words: }#1\par
}

\newtheorem{theorem}{Theorem}[section]
\newtheorem{lemma}[theorem]{Lemma}

\theoremstyle{definition}
\newtheorem{definition}[theorem]{Definition}

\theoremstyle{remark}

\title{TOPOLOGICAL PRESSURE OF DISCONTINUOUS POTENTIALS AND VARIATIONAL PRINCIPLE FOR FLOWS}
\author{
    Ruolan Xiong
    \and
    Yunhua Zhou}
\date{}

\begin{document}

\maketitle

\begin{abstract}
We introduce two types of topological pressures associated with a discontinuous potential $S=\left\{s_{t}\right\}_{t>0}$. Let $X$ be a compact metric space and $\Phi=\left\{\varphi_{t}\right\}_{t \in \mathbb{R}}$ be a continuous flow on $X$. The first type of topological pressure is defined for a family of measurable functions $S= \left\{s_{t}\right\}_{t>0}$, and its variational principle is established. The second type of topological pressure is defined in an invariant subset with a nested structure of subsets, and we demonstrate that it also satisfies a variational principle.
\end{abstract}

\maketitle

%\maketitle

\Keywords{topological pressure; discontinuous potentials; variational principle; flow
}

\section{Introduction}
Topological pressure, introduced in the 1970s, extends the concept of topological entropy--a measure of the complexity of motion in dynamical systems. The behavior of such systems varies under different potential functions. Thus, it is important to define a quantity that captures the relationship between the potential function and the uncertainty of the system. This quantity is known as topological pressure. 

The refined definition of topological pressure was initially formulated by David Ruelle, a prominent theoretical physicist and esteemed member of the French Academy of Sciences. Inspired by the earlier contributions of Sinai and Bowen, Ruelle introduced the concept of topological pressure in 1973, building upon the notion of ``pressure" from statistical mechanics. His definition was specifically applied to $\mathbb{Z}^{v}$--actions satisfying expansivity and specification on compact metric spaces. He established the following equation:
$$
P(\varphi)=\max _{\mu \in I}\{s(\mu)+\mu(\varphi)\},
$$commonly referred to as the variational principle ([21]). In 1984, Pesin and Pitskel extended this framework by defining the topological pressure for continuous mappings on noncompact sets ([17]). The concept of topological pressure for subadditive potential sequences had already emerged in the 1960s ([14]). In 1988, Falconer studied the topological pressure of sequences of subadditive potentials on mixed repellers, proving a variational principle under Lipschitz conditions and bounded variation ([10]). Subsequently, numerous studies have further investigated the topological pressure of sequences of subadditive potentials ([7, 9, 12, 23]). In 1996, Barreira relaxed the subadditivity condition and defined the topological pressure for sequences of potentials (not necessarily subadditive) on subsets of compact metric spaces ([5]). Over the years, the concept of topological pressure has evolved beyond continuous mappings in discrete dynamical systems. It has been extended to include the topological pressure and variational principles for families of continuous potentials $a=\left\{a_{t}\right\}_{t>0}$ ([4]). Similarly to measure-theoretic entropy, dynamical systems also incorporate measure-theoretic pressure, which originates from the concept of topological pressure and is closely related to invariant measures ([8, 13, 24]).

The above conclusions are primarily concerned with continuous potentials. However, when the potential does not satisfy continuity, does topological pressure still exist? In 2006, Mummert first defined the topological pressure for discontinuous potentials on a subset $\Lambda$ of a compact metric space ( $X, T$ ) ([16]). In Mummert's study, the set $\Lambda$ is represented as the union of a nested sequence of subsets: $\Lambda=\bigcup\limits_{l \geq 1} \Lambda_{l}$. The potential is required to be continuous only on the closure of each subset $\Lambda_{l}$, rather than on the entire set $\Lambda$. Based on this structure, the classical topological pressure $P_{\Lambda_{l}}(\varphi)$ can be defined on each subset $\Lambda_{l}$. By taking the supremum over $l$, the topological pressure on the subset $\Lambda$ is obtained as follows: $P_{\Lambda}(\varphi)=\sup \limits_{l \geq 1} P_{\Lambda_{l}}(\varphi)$. Subsequently, Ma Xianfeng et al. extended Mummert's results from a single potential to a sequence of subadditive potentials ([15]). In 2012, J. Barral and D. J. Feng investigated the topological pressure of upper semi-continuous subadditive potentials ([2]). In 2016, Feng and Huang introduced the concept of weighted topological pressure for upper semi-continuous entropy maps and established a corresponding variational principle ([11]). In 2017, Marc Rauch directly defined the topological pressure for measurable potentials on compact metric spaces using Caratheodory structure theory and proved the associated variational principle ([19]). In subsequent work, he further developed the topological pressure and variational principle for subadditive potential sequences ([20]). Additionally, topological pressure has also been defined for systems with discontinuous semi-flows ([1]).

There have been several results concerning the topological pressure of discontinuous potentials in discrete dynamical systems. However, in the context of continuous dynamical systems, this topic remains largely unexplored. This paper presents two definitions of the topological pressure for discontinuous potentials with respect to flows and establishes the corresponding variational principles.

First, the topological pressure of a family of measurable potentials on a subset is introduced. In a continuous dynamical system ( $X, \Phi$ ), the set of all $\Phi$-invariant probability measures on ( $X, \mathcal{B}$ ) is denoted by $\mathcal{M}_{\Phi}(X)$, and the set of ergodic probability measures is denoted by $\mathcal{E}_{\Phi}(X)$.

Given a nonempty subset $Z \subseteq X(Z \neq \emptyset)$, and real numbers $\epsilon>0$ and $\alpha \in \mathbb{R}$, we define the topological pressure of the family of potentials $S$ on the set $Z$ as
$$
P_{Z}(S)=\lim _{\epsilon \rightarrow 0} P_{Z}(S, \epsilon)=\lim _{\epsilon \rightarrow 0} \inf \{\alpha \in \mathbb{R}: M(Z, S, \alpha, \epsilon)=0\}
$$
where
$$
M(Z, S, \alpha, \epsilon)=\lim _{T \rightarrow+\infty} M(Z, S, \alpha, \epsilon, T),
$$
and
$$
M(Z, S, \alpha, \epsilon, T)=\inf _{\Gamma} \sum_{(x, t) \in \Gamma} \exp (S(x, t, \epsilon)-\alpha t) .
$$
Here, the infimum is taken over all countable open covers $\Gamma \subseteq X \times[T,+\infty)$ that cover $Z$, i.e., $Z \subseteq \bigcup\limits_{(x, t) \in \Gamma} B_{t}(x, \epsilon)$, and $S(x, t, \epsilon)=\sup \left\{s_{t}(y): y \in B_{t}(x, \epsilon)\right\}$.

In the second part of this paper, we demonstrate that this definition is well-defined, stemming from the Carath\'{e}odory dimension theory of Pesin ([18]). Finally, we present the first main result of this article.

\begin{theorem}
    Let $\Phi$ be a continuous flow on $X$. Define $\mathcal{G}(S, \lambda)$ as the set
$$
\mathcal{G}(S, \lambda):=\left\{\mu \in \mathcal{A}(\lambda) \cap \mathcal{E}_{\Phi}(X) \cap \mathcal{M}_{\Phi}(X): \lim _{t \rightarrow \infty} \frac{1}{t} s_{t}(x)=\lambda(\mu) \text { for } \mu \text {-almost every } x \in X\right\}
$$
For each subset $\mathcal{Y} \subseteq \mathcal{G}(S, \lambda)$, it holds that
$$
P_{A(S, \lambda, \mathcal{Y})}(S)=\sup \left\{h_{\mu}(\Phi)+\lambda(\mu): \mu \in \mathcal{Y}\right\}
$$
In particular, for each $\mu \in \mathcal{Y}$, there exists a Borel set $B_{\mu} \subseteq A(A, \lambda, \mathcal{Y})$ such that $\mu\left(B_{\mu}\right)=1$, and
$$
P_{B}(S)=\sup \left\{h_{\mu}(\Phi)+\lambda(\mu): \mu \in \mathcal{Y}\right\},
$$where $B=\bigcup\limits_{\mu\in \mathcal{Y}}B_\mu$.
\end{theorem}

The set $A(S, \lambda, \mathcal{Y})$ is a subset of $X$ which is defined in Section 3.

For the second definition, consider a $\Phi$-invariant subset $Z \subseteq X$. The set $Z$ is constructed as the union of a nested sequence of subsets $\left\{Z_{l}\right\}_{l \geq 1}$ such that $Z=\bigcup\limits_{l \geq 1} Z_{l}$ and $Z_{l} \subseteq Z_{l+1}$ for all $l \in \mathbb{N}$. We require that $S$ be continuous on the closure of each subset $Z_{l}$, but not necessarily on $Z$ (for example, systems with nonzero Lyapunov exponents (see 3.1 of [16])). The topological pressure of $S$ on $Z$ with respect to the flow $\Phi$ is defined as
$$
P_{Z}(S)=\sup _{l \geq 1} P_{Z_{l}}(S)
$$
where $P_{Z_{l}}(S)$ is the classical topological pressure of $S$ on $Z_{l}$. Since $S$ is continuous on the closure of every subset $Z_{l},  P_{Z_{l}}(S)$ is well-defined. We then obtain the second main result of this article.

\begin{theorem}
  Let $\Phi$ be a continuous flow on $X$ and $S=\left\{s_{t}\right\}_{t>0}$ be a family of functions with tempered variation such that
$$
\sup _{t \in[0, T]}\left\|s_{t}\right\|_{\infty}<+\infty \quad \text { for all } T>0 .
$$
Let $Z=\bigcup\limits_{l \geq 1} Z_{l} \subseteq X$ be a Borel $\Phi$-invariant set, and suppose $S$ is continuous with respect to the family of subsets $\left\{Z_{l}\right\}$. If there exists a continuous function $b: X \rightarrow \mathbb{R}$ such that
$$
s_{t+p}-s_{t} \circ \varphi(x) \rightarrow \int_{0}^{p}\left(b \circ \varphi_{u}\right)\, d u
$$
uniformly on $Z$ as $t \rightarrow+\infty$ for some $p>0$, then
$$
P_{\mathcal{L}(Z)}(S)=\sup \left\{h_{\mu}(\Phi)+\int_{Z} b\, d \mu: \mu \in \mathcal{M}_{Z}\right\}
$$  
\end{theorem}

Although the two topological pressures above are defined using different methods, both reflect the relationship between the potentials and the complexity of the system's motion on a subset of a continuous dynamical system.
\section{Preliminaries}
Let ( $X, d$ ) be a compact metric space and $\Phi=\left(\varphi_{t}\right)_{t \in \mathbb{R}}$ be a continuous flow on $X$.
A measure $\mu$ on ( $X, \mathcal{B}$ ) is said to be $\Phi$-invariant if $\mu\left(\varphi_{t}^{-1} B\right)=\mu(B)$ for every $B \in \mathcal{B}$ and $t \in \mathbb{R}$. A set $Z \subseteq X$ is called $\Phi$-invariant if $\varphi_{t}^{-1}(Z)=Z$ for all $t \in \mathbb{R}$. A measure $\mu$ on $(X, \mathcal{B})$ is said to be ergodic for the flow if every Borel $\Phi$-invariant set $B \in \mathcal{B}$ satisfies $\mu(B)=0$ or $\mu(B)=1$.

For each measurable subset $Z \subseteq X$, we define
$$
\mathcal{M}_{Z}=\left\{\mu \in \mathcal{M}_{\Phi}(X) \cap \mathcal{E}_{\Phi}(X): \mu(Z)=1\right\}.
$$
For $\mu \in \mathcal{M}_{\Phi}(X)$, the quantity $h_{\mu}(\Phi)$ is defined as $h_{\mu}\left(\varphi_{1}\right)$.

Given $x \in X$ and $t, \epsilon>0$, consider the set
$$
B_{t}(x, \epsilon)=\left\{y \in X: d\left(\varphi_{s}(x), \varphi_{s}(y)\right)<\epsilon, s \in[0,\, t]\right\}
$$
A collection $\mathcal{U} \subseteq 2^{X}$ is called a finite open cover of $X$ if $\#(\mathcal{U})<\infty$ and $X \subseteq \bigcup\limits_{U \in \mathcal{U}} U$, where each $U \in \mathcal{U}$ is open.

Given $x \in X$ and $t>0$, define the probability measures
$$
\delta_{x, t}=\frac{1}{t} \int_{0}^{t} \delta_{\varphi_{s}(x)} \,d s
$$
where $\delta_{\varphi_{s}(x)}$ is the probability measure concentrated on $\left\{\varphi_{s}(x)\right\}$. The set of all $\Phi$-invariant sublimits of $\left\{\delta_{x, t}\right\}_{t>0}$ in the weak*-topology is denoted by $V_{\Phi}(x)$. We know that $V_{\Phi}(x) \subseteq M_{\Phi}(X)$ and $V_{\Phi}(x) \neq \emptyset$ for every $x \in X$. Furthermore, $V_{\Phi}(x)$ is a compact metrizable space. Let
$$
Z_{\mu}:=\left\{x \in X: \lim _{t \rightarrow \infty} \delta_{x, t}=\mu\right\}
$$
Note that $\mu\left(Z_{\mu}\right)=1$ if $\mu$ is ergodic.

Given $\epsilon>0$, we say that a set $\Gamma \subseteq X \times \mathbb{R}_{0}^{+}$covers a subset $Z \subseteq X$ if $Z \subseteq \bigcup\limits_{(x, t) \in \Gamma} B_{t}(x, \epsilon)$. Let $\mathcal{G}_{Z}(\epsilon, T)$ be the set of all countable covers $\Gamma \subseteq X \times \mathbb{R}_{0}^{+}$ of $Z$ by Bowenballs $B_t(x,\,\epsilon)$. Let $S=\left\{s_{t}\right\}_{t>0}$ be a family of functions on $X$ satisfying
$$
\lim _{\epsilon \rightarrow 0} \varlimsup_{t \rightarrow \infty} \frac{\gamma_{t}(S, \epsilon)}{t}=0
$$
where
$$
\gamma_{t}(S, \epsilon)=\sup \left\{\left|s_{t}(y)-s_{t}(z)\right|: y, z \in B_{t}(x, \epsilon) \text { for some } x \in X\right\}
$$
which is called tempered variation. 

Fix $\emptyset \neq Z \subseteq X, \epsilon>0, \alpha \in \mathbb{R}$. Set
$$
M(Z, S, \alpha, \epsilon, T)=\inf _{\Gamma\in \mathcal{G}_{Z}(\epsilon, T)} \sum_{(x, t) \in \Gamma} \exp (S(x, t, \epsilon)-\alpha t),
$$
where the infimum is taken over all countable sets $\Gamma \subseteq X \times[T,+\infty)$ covering $Z$. Define
$$
M(Z, S, \alpha, \epsilon)=\lim _{T \rightarrow+\infty} M(Z, S, \alpha, \epsilon, T).
$$
It is clear that $M(Z, S, \alpha, \epsilon, T)$ increases as $T$ increases. Hence, we have
$$
M(Z, S, \alpha, \epsilon)=\sup _{T} M(Z, S, \alpha, \epsilon, T).
$$
\begin{lemma}
   Let $\beta \in \mathbb{R}$ and $Z \subseteq X$. If $M(Z, S, \beta, \epsilon)<\infty$, then $M(Z, S, \alpha, \epsilon)=0$ for all $\alpha>\beta$. 
\end{lemma}
\begin{proof}
  If $Z \neq \emptyset$, by definition, there exists a real number $M$ such that $M(Z, S, \beta, \epsilon)<M$, and consequently, $M(Z, S, \beta, \epsilon, T)<M$ for all $T>0$. Hence, we have
$$
\begin{aligned}
0 & \leq \inf _{\Gamma} \sum_{(x, t) \in \Gamma} \exp (S(x, t, \epsilon)-\alpha t) \\
& =\inf _{\Gamma} \sum_{(x, t) \in \Gamma} \exp (S(x, t, \epsilon)-\beta t) \exp (t(\beta-\alpha)) \\
& \leq M(\exp (\beta-\alpha))^{T}
\end{aligned}
$$

Since $\beta-\alpha \leq 0$, it follows that $M(\exp (\beta-\alpha))^{T} \rightarrow 0$ as $T \rightarrow \infty$. Therefore, $M(Z, S, \alpha, \epsilon)=0$ for all $\alpha>\beta$.  
\end{proof}

By Lemma 2.1, we can define the following quantity:
$$
P_{Z}(S, \epsilon)=\inf \{\alpha \in \mathbb{R}: M(Z, S, \alpha, \epsilon)=0\}
$$

For $0<\epsilon^{\prime}<\epsilon$ and $T \geq 0$, we have $B_{t}\left(x, \epsilon^{\prime}\right) \subseteq B_{t}(x, \epsilon)$ for all $t>T$. If $Z \subseteq \bigcup\limits_{(x, t) \in \Gamma} B_{t}\left(x, \epsilon^{\prime}\right)$, then $Z \subseteq \bigcup\limits_{(x, t) \in \Gamma} B_{t}(x, \epsilon)$. This shows that $\mathcal{G}_{Z}\left(\epsilon^{\prime}, T\right) \subseteq \mathcal{G}_{Z}(\epsilon, T)$.

Combining Theorem 1 in [4], we conclude that the following limit is well-defined:
$$
P_{Z}(S)=\lim _{\epsilon \rightarrow 0} P_{Z}(S, \epsilon)=\sup _{\epsilon>0} P_{Z}(S, \epsilon)
$$
\setcounter{theorem}{0}
\begin{definition}\label{def:2.1}
The quantity $P_{Z}(S)$ is called the topological pressure of $S$ on $Z$ with respect to $\Phi$.
\end{definition}
\addtocounter{theorem}{0} 
Note that $M(\emptyset, S, \alpha, \epsilon, T)=0$ for all $\epsilon>0, \alpha \in \mathbb{R}$, and $T \geq 0$. Therefore, $P_{\emptyset}(S)=-\infty$ for every potential $S$. Thus, $P_{Z}(S) \neq-\infty$ implies $Z \neq \emptyset$.

Next, we introduce some properties of topological pressure.
\begin{lemma}
  For $\epsilon>0, \alpha \in \mathbb{R}, T \geq 0$, and $Y \subseteq Z \subseteq X$, the following inequality holds: $$ M (Y, S, \alpha, \epsilon, T) \leq M(Z, S, \alpha, \epsilon, T)
$$
In particular, $P_{Y}(S, \epsilon) \leq P_{Z}(S, \epsilon)$ for all $\epsilon>0$, and thus $P_{Y}(S) \leq P_{Z}(S)$. 
\end{lemma}
\begin{proof}
    Let $\Gamma \subseteq X \times \mathbb{R}_{0}^{+}$be a cover of $Z$. Since $Y \subseteq Z$, the same cover $\Gamma$ also covers $Y$. Therefore, we have
$$
M(Y, S, \alpha, \epsilon, T) \leq M(Z, S, \alpha, \epsilon, T).
$$
Then it is easy to see that
$$
P_{Y}(S, \epsilon) \leq P_{Z}(S, \epsilon).
$$
Finally, by taking the limit as $\epsilon \rightarrow 0$, we obtain $P_{Y}(S) \leq P_{Z}(S)$.
\end{proof} 
\begin{lemma}
    Given a set $Z \subseteq X$, suppose $Z=\bigcup_{i \in I} Z_{i}$, where $I \subseteq \mathbb{N}$ and $Z_{i} \subseteq X$ for all $i \in I$. Then
$$
P_{Z}(S)=\sup _{i \in I} P_{Z_{i}}(S)
$$
\end{lemma}
\begin{proof}
   Since $Z_{i} \subseteq Z$, we have $P_{Z_{i}}(S) \leq P_{Z}(S)$ for each $i \in I$. Thus,
$$
P_{Z}(S) \geq \sup _{i \in I} P_{Z_{i}}(S) .
$$

Now, take $\alpha \in \mathbb{R}$ such that $\alpha>\sup _{i \in I} P_{Z_{i}}(S, \epsilon)$. Then $M\left(Z_{i}, S, \alpha, \epsilon\right)=0$ for each $i \in I$. Hence, for any $\delta>0$ and $T>0$, there exists $\Gamma_{i} \subseteq X \times[T,+\infty)$ covering $Z_{i}$ such that
$$
\sum_{(x, t) \in \Gamma_{i}} \exp (S(x, t, \epsilon)-\alpha t)<\frac{\delta}{2^{i}}
$$

Let $\Gamma=\bigcup\limits_{i \in I} \Gamma_{i}$. Then $\Gamma$ covers $Z$, and we have
$$
\sum_{(x, t) \in \Gamma} \exp (S(x, t, \epsilon)-\alpha t) \leq \sum_{i \in I} \sum_{(x, t) \in \Gamma_{i}} \exp (S(x, t, \epsilon)-\alpha t) \leq \sum_{i \in I} \frac{\delta}{2^{i}} \leq \delta .
$$

This implies $M(Z, S, \alpha, \epsilon) \leq \delta$, and since $\delta$ is arbitrary, $M(Z, S, \alpha, \epsilon)=0$. Therefore, $\alpha \geq P_{Z}(S, \epsilon)$. Letting $\alpha \rightarrow \sup _{i \in I} P_{Z_{i}}(S, \epsilon)$, we obtain $\sup _{i \in I} P_{Z_{i}}(S, \epsilon) \geq P_{Z}(S, \epsilon)$.

Finally, letting $\epsilon \rightarrow 0$, we conclude that $\sup _{i \in I} P_{Z_{i}}(S) \geq P_{Z}(S)$. Combining the inequalities, we have
$$
P_{Z}(S)=\sup _{i \in I} P_{Z_{i}}(S)
$$ 
\end{proof}
\section{Variational principle for family of measurable potentials}
In this section, we introduce the variational principle for the topological pressure defined above.

For a mapping $\lambda: \mathcal{M}_{\Phi}(X) \rightarrow[-\infty,+\infty]$, the set
$$
\mathcal{A}(\lambda):=\left\{\mu \in \mathcal{M}_{\Phi}(X): h_{\mu}(\Phi)<\infty \text { or } \lambda(\mu)>-\infty\right\}
$$
is called the set of all allowed $\Phi$-invariant measures with respect to $\lambda$. For every $\mu \in \mathcal{A}(\lambda)$, the quantity $h_{\mu}(\Phi)+\lambda(\mu)$ is well-defined.

For $\mu \in \mathcal{M}_{\Phi}(X)$, a point $x \in X$ is called an allowed point with respect to $\lambda, a$, and $\mu$ if $\mu \in V_{\Phi}(x)$ and
\begin{equation*}
\limsup _{l \rightarrow \infty} \frac{1}{t_{l}} s_{t_{l}}(x) \leq \lambda(\mu) \tag{3.1}
\end{equation*}
for every sub-family $\left(t_{l}\right)_{l>0}$ satisfying $\delta_{x, t_{l}} \rightarrow \mu$ as $l \rightarrow \infty$. The set of all such points $x$ is denoted by $A(S, \lambda, \mu)$.

For a subset $\mathcal{Y} \subseteq \mathcal{M}_{\Phi}(X)$, define
$$
A(S, \lambda, \mathcal{Y}):=\bigcup_{\mu \in \mathcal{Y}} A(S, \lambda, \mu), \quad V(S, \lambda, \mathcal{Y}):=\bigcap_{\mu \in \mathcal{Y}} A(S, \lambda, \mu)
$$

In fact, these two sets can often be expressed more concisely under certain conditions. The following Lemma 3.1 illustrates this point.

\begin{lemma}
    Let $f: X \rightarrow[-\infty,+\infty]$ be an upper semi-continuous function, $\lambda(\mu):= \int_{X} f d \mu$, and $s_{t}:=\int_{0}^{t} f \circ \varphi_{u} d u$. Then, the following hold:
$$
A(S, \lambda, \mathcal{Y})=\left\{x \in X: V_{\Phi}(x) \cap \mathcal{Y} \neq \emptyset\right\}, \quad V(S, \lambda, \mathcal{Y})=\left\{x \in X: \mathcal{Y} \subseteq V_{\Phi}(x)\right\}
$$
\end{lemma}
\begin{proof}
   Fix $x \in X$ and corresponding $\mu \in V_{\Phi}(x) \cap \mathcal{Y}$. This implies $x \in\{x \in X$ : $\left.V_{\Phi}(x) \cap \mathcal{Y} \neq \emptyset\right\}$. For a sub-family $\left(t_{l}\right)_{l>0}$ such that $\lim _{l \rightarrow \infty} \delta_{x, t_{l}}=\mu$, we have, since $f$ is upper semi-continuous, by Lemma A.2(d) in [20],
$$
\limsup\limits _{l \rightarrow \infty} \frac{1}{t_{l}} s_{t_{l}}(x)=\limsup \limits_{l \rightarrow \infty} \frac{1}{t_{l}} \int_{0}^{t_{l}} f\left(\varphi_{u}(x)\right)\, d u=\limsup\limits _{l \rightarrow \infty} \int_{X} f\, d \delta_{x, t_{l}} \leq \int_{X} f \,d \mu.
$$
Thus, $x \in A(S, \lambda, \mu)$, and consequently, $\left\{x \in X: V_{\Phi}(x) \cap \mathcal{Y} \neq \emptyset\right\} \subseteq A(S, \lambda, \mathcal{Y})$.

Conversely, if $x \in A(S, \lambda, \mathcal{Y})$, there exists $\mu \in \mathcal{Y}$ such that $x \in A(S, \lambda, \mu)$. By definition, this implies $V_{\Phi}(x) \cap \mathcal{Y} \neq \emptyset$. Therefore, $A(S, \lambda, \mathcal{Y}) \subseteq\left\{x \in X: V_{\Phi}(x) \cap \mathcal{Y} \neq \emptyset\right\}$, and we conclude that
$$
A(S, \lambda, \mathcal{Y})=\left\{x \in X: V_{\Phi}(x) \cap \mathcal{Y} \neq \emptyset\right\}.
$$
For the second equation, let $x \in V(S, \lambda, \mathcal{Y})$. Then, for every $\mu \in \mathcal{Y}$, since $x \in A(S, \lambda, \mu)$, we have $\mu \in V_{\Phi}(x)$. This implies $\mathcal{Y} \subseteq V_{\Phi}(x)$, and hence $V(S, \lambda, \mathcal{Y}) \subseteq\left\{x \in X: \mathcal{Y} \subseteq V_{\Phi}(x)\right\}$.

Conversely, if $x \in\left\{x \in X: \mathcal{Y} \subseteq V_{\Phi}(x)\right\}$, then for every $\mu \in \mathcal{Y}$, we have $\mu \in V_{\Phi}(x)$. Hence, there exists a subsequence $\left(t_{l}\right)_{l>0}$ such that $\lim\limits _{l \rightarrow \infty} \delta_{x, t_{l}}=\mu$. Similarly,
$$
\limsup _{l \rightarrow \infty} \frac{1}{t_{l}} s_{t_{l}}(x) \leq \int_{X} f \,d \mu=\lambda(\mu)
$$
Thus, $x \in V(S, \lambda, \mu)$ for every $\mu \in \mathcal{Y}$, which implies $x \in V(S, \lambda, \mathcal{Y})$. Therefore, $V(S, \lambda, \mathcal{Y}) \supseteq \left\{x \in X: \mathcal{Y} \subseteq V_{\Phi}(x)\right\}$.

Combining the above, we conclude that
$$
V(S, \lambda, \mathcal{Y})=\left\{x \in X: \mathcal{Y} \subseteq V_{\Phi}(x)\right\}
$$ 
\end{proof}

Before presenting the variational principle, we first introduce two auxiliary lemmas.
\begin{lemma}
    ([6]) Let $E$ be a finite set, given $q \in \mathbb{N}$ and $\alpha=\left(a_{1}, a_{2}, \ldots a_{q}\right) \in E^{q}$. Define the probability measure $\nu_{\alpha}$ as
$$
\nu_{\alpha}(e)=\frac{1}{q} \#\left\{j \in \mathbb{N}: a_{j}=e\right\}
$$
for every $e \in E$, and set the entropy of $\alpha$ to
$$
H(\alpha):=-\sum_{e \in E} \nu_{\alpha}(e) \log \nu_{\alpha}(e)
$$
Then for $h \geq 0$, one has
$$
\limsup _{q \rightarrow \infty} \frac{1}{q} \log \#\left\{\alpha \in E^{q}: H(\alpha) \leq h\right\} \leq h.
$$
\end{lemma}
\begin{lemma}
   Let $x \in A(S, \lambda, \mathcal{Y})$ and $\mu \in V_{\Phi}(x)$ such that $x$ is allowed with respect to $\lambda, S$, and $\mu$. Let $\delta>0$, and let $\Gamma \subseteq X \times\{1\}$ be a finite cover of $X$. For the open cover $\mathcal{U}=\left\{V_{1}, V_{2}, \ldots, V_{r}\right\}$ of $X$, where $V_{j}=B_{1}\left(x_{j}, \frac{\epsilon}{2}\right)$ with $\left(x_{j}, 1\right) \in \Gamma$, there exist integers $m, p \in \mathbb{N}$, with $p$ arbitrarily large, and a sequence $U=V_{i_{1}} V_{i_{2}} \ldots V_{i_{p}}$ such that:\\
(a) $x \in \bigcap_{r=1}^{p} \varphi_{-r+1}\left(V_{i_{r}}\right)$.\\
(b) There exists a subset $V \in\left(\mathcal{U}^{m}\right)^{q}$ of $U$, of length $p \geq q m \geq p-m$, satisfying the inequality
$$
H(V) \leq m\left(h_{\mu}(\Phi)+\delta\right).
$$
(c) $a_{p}(x) \leq p(\lambda(\mu)+\delta)$. 
\end{lemma}
\begin{proof}
   The statements (a) and (b) are given in [4]. We will now prove statement (c) in a manner similar to Lemma 4.7 (3) in [20].

There exists an increasing sub-family $\left(t_{j}^{\prime}\right)_{j \in \mathbb{R}}$ and a corresponding sequence $\left\{V_{i_{j}}\right\}_{j=1}^{p}$ satisfying (a) and (b), with $\delta_{x, t_{j}^{\prime}} \rightarrow \mu$ as $j \rightarrow \infty$. By inequality (3.1), we have
$$
\limsup _{j \rightarrow \infty} \frac{1}{t_{j}^{\prime}} s_{t_{j}^{\prime}}(x) \leq \lambda(\mu).
$$
Thus, there exists an index $j_{0} \geq 0$ such that
$$
s_{t_{j}^{\prime}}(x) \leq t_{j}^{\prime}(\lambda(\mu)+\delta) \quad \text { for all } j \geq j_{0}
.$$
For each $j \geq 0$, let $p:=\left[t_{j+j_{0}}^{\prime}\right]+1$. Together with the corresponding sequence $V_{i_{j+j_{0}}}$, this satisfies (c). 
\end{proof}
\newcounter{tempcounter}
\setcounter{tempcounter}{\value{theorem}}
\setcounter{theorem}{0}
\begin{theorem}\label{thm:3.1}
    For a mapping $\lambda: \mathcal{M}_{\Phi}(X) \rightarrow[-\infty,+\infty]$ and a subset $\mathcal{Y} \subseteq \mathcal{A}(\lambda)$, the following inequalities hold:
\begin{align*}
& P_{A(S, \lambda, \mathcal{Y})}(S) \leq \sup \left\{h_{\mu}(\Phi)+\lambda(\mu): \mu \in \mathcal{Y}\right\},  \tag{3.2}\\
& P_{V(S, \lambda, \mathcal{Y})}(S) \leq \inf \left\{h_{\mu}(\Phi)+\lambda(\mu): \mu \in \mathcal{Y}\right\}. \tag{3.3}
\end{align*}
\end{theorem}
\begin{proof}
 Let $\sup \emptyset:=-\infty$. For each $\mathcal{Y} \subseteq \mathcal{A}(\lambda)$, for simplicity, let $P_{\mathcal{Y}}(\lambda):=\sup \left\{h_{\mu}(\Phi)+\right. \lambda(\mu): \mu \in \mathcal{Y}\}$ and $Q_{\mathcal{Y}}(\lambda):=\inf \left\{h_{\mu}(\Phi)+\lambda(\mu): \mu \in \mathcal{Y}\right\}$. If $P_{\mathcal{Y}}(\lambda)=\infty$, the conclusion is clear. Now assume that $P_{\mathcal{Y}}(\lambda)<\infty$. This implies $0 \leq h_{\mu}(\Phi)<\infty$ and $\lambda(\mu)<\infty$ for each $\mu \in \mathcal{Y}$. Denote $\mathcal{Y}_{-\infty}:=\{\mu \in \mathcal{Y}: \lambda=-\infty\}$ and $\mathcal{Y}^{\prime}:=\{\mu \in \mathcal{Y}: \lambda>-\infty\}$, then $\mathcal{Y}=\mathcal{Y}_{-\infty} \cup \mathcal{Y}^{\prime}$. Hence, $A(S, \lambda, \mathcal{Y})=A\left(S, \lambda, \mathcal{Y}_{-\infty}\right) \cup A\left(S, \lambda, \mathcal{Y}^{\prime}\right)$. By Lemma 2.3, we obtain
   \begin{equation*}
P_{A(S, \lambda, \mathcal{Y})}(S)=\max \left\{P_{A\left(S, \lambda, \mathcal{Y}_{-\infty}\right)}(S), P_{A\left(S, \lambda, \mathcal{Y}^{\prime}\right)}(S)\right\}. \tag{3.4}
\end{equation*}

For the subsequent proof, we assume that (3.2) holds for each $\mathcal{Y} \subseteq \mathcal{A}(\lambda)$ such that $\lambda(\mu) \in (-\infty, \infty)$ for all $\mu \in \mathcal{Y}$. Define a sequence of mappings $\lambda_{N}(\mu):=-h_{\mu}(\Phi)-N$ for $N \in \mathbb{N}$ and $\mu \in \mathcal{M}_{\Phi}(X)$. For $\mu \in \mathcal{Y}_{-\infty}$, one has $-\infty<\lambda_{N}(\mu)<\infty$ for all $N \in \mathbb{N}$ since $h_{\mu}(\Phi)<\infty$, and thus $A\left(S, \lambda, \mathcal{Y}_{-\infty}\right) \subseteq A\left(S, \lambda_{N}, \mathcal{Y}_{-\infty}\right)$. If $\mathcal{Y}_{-\infty}=\emptyset$, we have
$$
P_{A\left(S, \lambda, \mathcal{Y}_{-\infty}\right)}(S)=P_{\mathcal{Y}_{-\infty}}\left(\lambda_{N}\right)=-\infty
$$
If $\mathcal{Y}_{-\infty} \neq \emptyset$, by (3.2) we obtain
$$
\begin{aligned}
P_{A\left(S, \lambda, \mathcal{Y}_{-\infty}\right)}(S) & \leq P_{A\left(S, \lambda_{N}, \mathcal{Y}_{-\infty}\right)}(S) \leq P_{\mathcal{Y}_{-\infty}}\left(\lambda_{N}\right) \\
& =\sup \left\{h_{\mu}(\Phi)-h_{\mu}(\Phi)-N: \mu \in \mathcal{Y}_{-\infty}\right\} \\
& \leq-N.
\end{aligned}
$$
Letting $N \rightarrow \infty$, we have $P_{A\left(S, \lambda, \mathcal{Y}_{-\infty}\right)}(S)=-\infty$. By (3.4), we obtain
$$
P_{A(S, \lambda, \mathcal{Y})}(S)=P_{A\left(S, \lambda, \mathcal{Y}^{\prime}\right)}(S) \leq P_{\mathcal{Y}^{\prime}}(\lambda) \leq P_{\mathcal{Y}}(\lambda).
$$
Next, we prove that (3.2) holds for all subsets $\mathcal{Y} \subseteq \mathcal{A}(\lambda)$ satisfying $h_{\mu}(\Phi)<\infty$ and $-\infty< \lambda(\mu)<\infty$ for all $\mu \in \mathcal{Y}$. Choose such a set $\mathcal{Y}$. If $P_{\mathcal{Y}}(\lambda)=-\infty$, then $\mathcal{Y}=\emptyset$. In fact, if there exists $\mu \in \mathcal{Y}$ such that $P_{\mathcal{Y}}(\lambda)=-\infty$, then
$$
P_{\mathcal{Y}}(\lambda) \geq h_{\mu}(\Phi)+\lambda(\mu)>-\infty
$$
since $h_{\mu}(\Phi)<\infty$ and $-\infty<\lambda(\mu)<\infty$, which is a contradiction. Hence, $\mathcal{Y}=\emptyset$, and thus $A(a, \lambda, \mathcal{Y})=\emptyset$. Therefore, $P_{A(S, \lambda, \mathcal{Y})}(S)=P_{\mathcal{Y}}(\lambda)=-\infty$.

Consequently, we may assume
\begin{equation*}
-\infty<P_{\mathcal{Y}}(\lambda)<\infty .\tag{3.5}
\end{equation*}

Let $\delta>0$. For an open cover $\mathcal{U}$ of $X$ as in Lemma 3.3, pick $\mu_{x} \in \mathcal{Y}$ for each $x \in A(S, \lambda, \mathcal{Y})$ such that $x \in A\left(S, \lambda, \mu_{x}\right)$. We construct a sequence of real numbers $\left\{u_{i}\right\}_{i \geq 1}$ such that $\left\{u_{i}\right\}_{i \geq 1}$ is $\delta$-dense in $\mathbb{R}$, i.e., for any $y \in \mathbb{R}$, there exists an $i$ such that $\left|u_{i}-y\right|<\delta$. For $m, i \geq 1$, define the subset $A_{m, i} \subseteq A(S, \lambda, \mathcal{Y})$ as follows: a point $x \in A(S, \lambda, \mathcal{Y})$ is in $A_{m, i}$ if $x, m, \mu_{x}, \delta$, and $\mathcal{U}$ satisfy conditions (a), (b), and (c) in Lemma 3.3, and $\lambda\left(\mu_{x}\right) \in\left[u_{i}-\delta, u_{i}+\delta\right]$.

For each $x \in A(S, \lambda, \mathcal{Y}), \lambda\left(\mu_{x}\right) \in \mathbb{R}$ by (3.5), and hence there exists an $i$ such that $\lambda\left(\mu_{x}\right) \in \left[u_{i}-\delta, u_{i}+\delta\right]$. Combining this with Lemma 3.3, there exists a corresponding $m$ such that $x \in A_{m, i}$. Thus, the sets $A_{m, i}$ cover $A(S, \lambda, \mathcal{Y})$, i.e.,
$$
A(S, \lambda, \mathcal{Y})=\bigcup_{m \in \mathbb{N}} \bigcup_{i \in \mathbb{N}} A_{m, i}.
$$

For $N \geq 1$, define
$$
D_{N}:=\# \bigcup_{x \in A_{m, i}}\left\{U \in \mathcal{U}^{N}: U \text { satisfies }(a),(b),(c)\right\}.
$$
For $A_{m, i}$ and $q \geq 1$, we denote
$$
R_{q}:=\left\{V \in\left(\mathcal{U}^{m}\right)^{q}: H(V) \leq m\left(P_{\mathcal{Y}}(\lambda)-u_{i}+2 \delta\right)\right\} .
$$
Choose some $x \in A_{m, i}$. By Lemma 3.3, there exists $N \geq 1$ and a corresponding $q \leq \frac{N}{m}$, $U \in \mathcal{U}^{N}$, and $V \in\left(\mathcal{U}^{m}\right)^{q}$ such that
$$
0 \leq \frac{1}{m} H(V) \leq h_{\mu_{x}}(\Phi)+\delta \leq h_{\mu_{x}}(\Phi)+\lambda\left(\mu_{x}\right)-u_{i}+2 \delta \leq P_{\mathcal{Y}}(\lambda)-u_{i}+2 \delta
$$
Hence, $V \in R_{q}$, and since $m\left(P_{\mathcal{Y}}(\lambda)-u_{i}+2 \delta\right) \geq 0$, by Lemma 3.2 we obtain
$$
\limsup _{q \rightarrow \infty} \frac{1}{q} \log \#\left\{\alpha \in E^{q}: H(\alpha) \leq P_{\mathcal{Y}}(\lambda)-u_{i}+2 \delta\right\} \leq P_{\mathcal{Y}}(\lambda)-u_{i}+2 \delta
$$

Thus, there exists $q_{0} \in \mathbb{N}$ such that for all $q \geq q_{0}$,
$$
\frac{1}{q} \log \# R_{q} \leq m\left(P_{\mathcal{Y}}(\lambda)-u_{i}+3 \delta\right)
$$
Fix $N \geq N_{0}:=q_{0} m$. We now estimate $D_{N}$. For $q \geq q_{0}$ and each $V \in\left(\mathcal{U}^{m}\right)^{q}, V \in R_{q}$ if $\mathcal{U}$ satisfies the properties of $D_{N}$. The difference between $N$ and $m q$ is at most $m$ by condition (b) in Lemma 3.3. Hence,
$$
D_{N} \leq(\# \mathcal{U})^{m}\left(\# R_{q}\right) \leq(\# \mathcal{U})^{m} \exp \left(q m\left(P_ \mathcal{Y}(\lambda)-u_{i}+3 \delta\right)\right)
.$$
Thus, since $N=q m+r$ for some $0 \leq r \leq m$, we have
$$
\limsup _{N \rightarrow \infty} \frac{1}{N} \log D_{N} \leq \limsup _{N \rightarrow \infty} \frac{1}{N}\left(m \log (\# \mathcal{U})+q m\left(P_{\mathcal{Y}}(\lambda)-u_{i}+3 \delta\right)\right) \leq P_{\mathcal{Y}}(\lambda)-u_{i}+3 \delta .
$$
Therefore, there exists $N_{1} \geq N_{0}$ such that for all $N \geq N_{1}$, 
$$
D_{N} \leq \exp \left(N\left(P_{\mathcal{Y}}(\lambda)-u_{i}+4 \delta\right)\right).
$$
For each \( l \geq N_1 \), define the set
\[
\Gamma_l := \left\{ U \in \bigcup\nolimits_{N \geq l} \mathcal{U}^N : U \text{ satisfies the properties of Lemma~3.3 for some } x \in A_{m,i} \right\}.
\]
Let \(\Gamma'_l\) denote the set of all corresponding pairs \((x, t)\), where \(t = |U|\) is the length of the sequence \(U\). Then \(\Gamma'_l\) forms a cover of $A_{m,i}$. Note that, by Lemma 3.3, for each $U \in \Gamma_{l}$, we have
$$
s_{N^{\prime}}(x) \leq N^{\prime}\left(u_{i}+\delta\right)
,$$
where $N^{\prime}=\lvert U\rvert$ denotes the length of the sequence $U$.

Hence, for $\alpha \in \mathbb{R}$ and $l \geq N_{1}$, we obtain
$$
\begin{aligned}
M\left(Z_{m, i}, S, \alpha, \epsilon, l\right) & \leq \sum_{(x, t) \in \Gamma_{l}^{\prime}} \exp \left(-\alpha t+t\left(u_{i}+\delta+\frac{\gamma_{t}(S, \epsilon)}{t}\right)\right) \\
& \leq \sum_{N=l}^{\infty} D_{N} \exp \left(-\alpha N+N\left(u_{i}+\delta+\varlimsup_{t \rightarrow \infty} \frac{\gamma_{t}(S, \epsilon)}{t}\right)\right) \\
& \leq \sum_{N=l}^{\infty}\left(\exp \left(-\alpha+P_{\mathcal{Y}}(\lambda)+5 \delta+\varlimsup_{t \rightarrow \infty} \frac{\gamma_{t}(S, \epsilon)}{t}\right)\right)^{N}
\end{aligned}
$$
Since $\epsilon$ will eventually become sufficiently small, we may assume without loss of generality that $\alpha>P_{\mathcal{Y}}(\lambda)+5 \delta$. Let $K:=\exp \left(-\alpha+P_{\mathcal{Y}}(\lambda)+5 \delta\right)$, then $K<1$, and consequently,
$$
M\left(A_{m, i}, S, \alpha, \epsilon\right) \leq \limsup _{l \rightarrow \infty} \sum_{N=l}^{\infty} K^{N}=0 .
$$
This implies
$$
P_{A_{m, i}}(S, \epsilon) \leq P_{\mathcal{Y}}(\lambda)+5 \delta+\varlimsup_{t \rightarrow \infty} \frac{\gamma_{t}(S, \epsilon)}{t}
$$
for fixed $A_{m, i}$. Taking the supremum over all $m, i$ yields
$$
P_{A(S, \lambda, \mathcal{Y})}(S, \epsilon)=\sup _{m, i} P_{A_{m, i}}(S, \epsilon) \leq P_{\mathcal{Y}}(\lambda)+5 \delta+\varlimsup_{t \rightarrow \infty} \frac{\gamma_{t}(S, \epsilon)}{t} .
$$
Note that $\lim _{\epsilon \rightarrow 0} \varlimsup_{\overline{\lim }}^{t \rightarrow \infty} \frac{\gamma_{t}(S, \epsilon)}{t}=0$ with tempered variation, so letting $\delta \rightarrow 0$ and $\epsilon \rightarrow 0$ results in
$$
P_{A(S, \lambda, \mathcal{Y})}(S) \leq P_{\mathcal{Y}}(\lambda) .
$$

For the statement (3.3), fix $\mu\in \mathcal{Y}$, $V(S, \lambda, \mathcal{Y}) \subseteq A(S, \lambda, \mu)$, so by Lemma 2.2 and (3.2), we have
$$
P_{V(S, \lambda, \mathcal{Y})}(S) \leq P_{A(S, \lambda, \mu)}(S) \leq h_{\mu}(\Phi)+\lambda(\mu) .
$$
Taking the infimum over all $\mu \in \mathcal{Y}$, we obtain the result.
\end{proof}
\setcounter{theorem}{\value{tempcounter}}
\begin{lemma}
    ([22]) Let ( $X, d$ ) be a compact metric space and $\Phi$ a continuous flow on $X$. For any $\mu \in \mathcal{E}_{\Phi}(X)$, define $h_{\mu}(x, \epsilon, t):=-\frac{1}{t} \log \mu\left(B_{t}(x, \epsilon)\right)$. Then one has
$$
\lim _{\epsilon \rightarrow 0} \liminf _{t \rightarrow \infty} h_{\mu}(x, \epsilon, t)=\lim _{\epsilon \rightarrow 0} \limsup _{t \rightarrow \infty} h_{\mu}(x, \epsilon, t)=h_{\mu}(\Phi)
$$
for $\mu$-almost $x \in X$.
\end{lemma}
 
By Lemma 3.4 and Theorem 4.8 of[20] we can obtain following Lemma.
\begin{lemma}
 Let $\Phi$ be a continuous flow on $X$. Fix $\mu \in \mathcal{E}_{\Phi}(X)$, and let $\left(S_{t}\right)_{t>0}$ be a Borel measurable potential on ( $X, \Phi$ ). Suppose there exists a constant $b \in[-\infty,+\infty]$ and a Borel set $B \subseteq X$ with $\mu(B)>0$, such that for each $x \in B$,
$$
\liminf _{t \rightarrow \infty} \frac{1}{t} s_{t}(x) \geq b .
$$
Then, if $h_{\mu}(\Phi)+b$ is well-defined, it holds that
$$
P_{B}(S) \geq h_{\mu}(\Phi)+b
.$$   
\end{lemma}

Next we prove the Theorem 1.1.

\begin{proof}
    By Theorem 3.1, $P_{A(S, \lambda, \mathcal{Y})}(S) \leq P_{\mathcal{Y}}(\lambda)$. For each $x \in A(S, \lambda, \mathcal{Y})$, since $\lim \limits_{t \rightarrow \infty} \frac{1}{t} s_{t}(x)= \lambda(\mu)$, we have $\liminf\limits _{t \rightarrow \infty} \frac{1}{t} s_{t}(x) \geq \lambda(\mu)$ for each $x \in A(S, \lambda, \mathcal{Y})$. Applying Lemma 3.5, we obtain $P_{A(S, \lambda, \mathcal{Y})}(S) \geq P_{\mathcal{Y}}(\lambda)$, and thus $P_{A(S, \lambda, \mathcal{Y})}(S)=P_{\mathcal{Y}}(\lambda)$.

For each $\mu \in \mathcal{Y}$, since $\mu$ is ergodic, there exists a Borel set $B_{\mu} \subseteq Z_{\mu}=\left\{x \in X: \delta_{t, x} \rightarrow \mu\right\}$ such that $\mu\left(B_{\mu}\right)=1$ and $\lim\limits_{t \rightarrow \infty} \frac{1}{t} s_{t}(x)=\lambda(\mu)$ for all $x \in B_{\mu}$. Therefore, $B_{\mu} \subseteq A(S, \lambda, \mu) \subseteq A(S, \lambda, \mathcal{Y})$, and by Lemma 3.5, this implies that

$$
h_{\mu}(\Phi)+\lambda(\mu)=P_{\{\mu\}}(\lambda) \leq P_{B_{\mu}}(S) \leq P_{A(S, \lambda, \mathcal{Y})}(S)
.$$
Taking the supremum on the left-hand side yields the result.
\end{proof}
\section{Topological pressure based on nested sets and variational principle}
We introduce a different definition of topological pressure for discontinuous potential\\$\mathrm{}\left(s_{t}\right)_{t>0}$ in this section. Before that, let's review the definition of topological pressure for continuous potentials on a compact metric space ( $X, d$ ).

Let $\Phi$ be a continuous flow on ( $X, d$ ), and let $S=\left(s_{t}\right)_{t>0}$ be a family of continuous functions $s_{t}: X \rightarrow \mathbb{R}$ with tempered variation. Given $\epsilon>0$, for each $Z \subseteq X$ and $\alpha \in \mathbb{R}$, we define
\begin{equation*}
M(Z, S, \alpha, \epsilon)=\lim _{T \rightarrow \infty} \inf _{\Gamma} \sum_{(x, t) \in \Gamma} \exp (S(x, t, \epsilon)-\alpha t), \tag{4.1}
\end{equation*}
where the infimum is taken over all countable sets $\Gamma \subseteq X \times[T,+\infty)$ covering $Z$. We then define
$$
P_{Z}(S, \epsilon)=\inf \{\alpha \in \mathbb{R}: M(Z, S, \alpha, \epsilon)=0\}
.$$

By Theorem 1 in [4], the quantity
$$
P_{Z}(S)=\lim _{\epsilon \rightarrow 0} P_{Z}(S, \epsilon)
$$
is well-defined and is called the non-additive topological pressure of the family $S$ on $Z$. Clearly, if $Z_{1} \subseteq Z_{2}$, then $P_{Z_{1}}(S) \leq P_{Z_{2}}(S)$.

Next, we introduce the topological pressure of discontinuous potentials. Suppose that $Z \subseteq X$ is $\Phi$-invariant, with a nested family of subsets $\left\{Z_{l}\right\}_{l \geq 1}$.
\begin{definition}
 A family of measurable functions $S=\left(s_{t}\right)_{t>0}: X \rightarrow \mathbb{R}$ is said to be continuous with respect to the family of subsets $\left\{Z_{l}\right\}_{l \geq 1}$ if $s_{t}$ is continuous on the closure of each $Z_{l}$ for all $t>0$.   
\end{definition}

Note that the potential $S$ is not necessarily continuous on $Z$. We define the topological pressure of $S$ on $Z$, with respect to $\Phi$, as
$$
P_{Z}(S)=\sup _{l \geq 1} P_{Z_{l}}(S)
$$
where $P_{Z_{l}}(S)$ is the topological pressure of $S$ on $Z_{l}$, as defined above. Although the topological pressure $P_{Z}(S)$ and the nested family of subsets $\left\{Z_{l}\right\}_{l \geq 1}$ are closely related by definition, choosing different families $\left\{Z_{l}\right\}_{l \geq 1}$ within the same $Z$ does not affect the value of $P_{Z}(S)$. Using Theorem 2 in [16], we can obtain the following Lemma.
\setcounter{theorem}{0}
\begin{lemma}
    Let $\left\{Z_{l}^{1}\right\}$ and $\left\{Z_{l}^{2}\right\}$ be nested families of subsets that exhaust the $\Phi$-invariant set $Z$, and let $\left\{s_{t}: X \rightarrow \mathbb{R}\right\}_{t>0}$ be continuous with respect to both $\left\{Z_{l}^{1}\right\}$ and $\left\{Z_{l}^{2}\right\}$. Then
$$
P_{Z}(S)=\sup _{l \geq 1} P_{Z_{l}^{1}}(S)=\sup _{l \geq 1} P_{Z_{l}^{2}}(S)
$$
\end{lemma}

Given a Borel $\Phi$-invariant set $Z \subseteq X, x \in X$, and $t>0$, consider the Borel $\Phi$-invariant set
$$
\mathcal{L}(Z)=\left\{x \in Z: V_{\Phi}(x) \cap \mathcal{M}_{Z} \neq \emptyset\right\}
$$

We now state the following Lemma, given by Barreira in [4], which will be used to prove Theorem 1.2.
\begin{lemma}
    Given $x \in \mathcal{L}(Z) \cap Z_{l}$ and $\mu \in V_{\Phi}(x) \cap \mathcal{M}_{Z}$, let $\Gamma \subseteq X \times\{1\}$ be a finite cover of $X$ for the open cover $\mathcal{U}=\left\{V_{1}, V_{2}, \ldots, V_{r}\right\}$ of $X$, where $V_{j}=B_{1}\left(x_{j}, \frac{\epsilon}{2}\right)$ with $\left(x_{j}, 1\right) \in \Gamma$. Then, there exist $m, p \in \mathbb{N}$ with $p$ arbitrarily large, and a sequence $U=V_{i_{1}} V_{i_{2}} \ldots V_{i_{p}}$ such that\\
(a) $x \in \bigcap_{r=1}^{p} \varphi_{-r+1} V_{i_{r}}$ and $s_{p}(x) \leq p\left(\int_{Z} b d \mu+3 \delta\right)$ where $b$ is the continuous function from the hypothesis of Theorem 1.2.\\
(b) There exists a subset $V \in\left(\mathcal{U}^{m}\right)^{k}$ of $U$ of length $k m \geq p-m$ satisfying the inequality
$$
H(V) \leq m\left(h_{\mu}(\Phi)+\delta\right)
$$
\end{lemma}

We now begin the proof of Theorem 1.2.

\begin{proof}
    Firstly, we show that $P_{\mathcal{L}(Z)}(S) \leq \sup \left\{h_{\mu}(\Phi)+\int_{Z} b d \mu: \mu \in \mathcal{M}_{Z}\right\}$. Since $\mathcal{L}(Z)= \bigcup_{l \geq 1}\left(\mathcal{L}(Z) \cap Z_{l}\right)$, we have that
$$
P_{\mathcal{L}(Z)}(S)=\sup _{l \geq 1} P_{\mathcal{L}(Z) \cap Z_{l}}(S)
.$$
We now show that for every $l \geq 1, P_{\mathcal{L}(Z) \cap Z_{l}}(S) \leq \sup \left\{h_{\mu}(\Phi)+\int_{Z} b d \mu: \mu \in \mathcal{M}_{Z}\right\}$. 

Given a sequence $\{u_1,\,u_2,\,\cdots,\,u_r\}$ that is $\delta-$ dense in $[-M,\,M]$, where $M$ is a real number with $M>\left|b\right|$. Let $B_{m, i}$ be the set of $x \in \mathcal{L}(Z) \cap Z_{l}$ satisfying conditions (a) and (b) in Lemma 4.2 for the corresponding $\mu \in V_{\Phi}(x) \cap \mathcal{M}_{Z}$ and $\int_{Z} b \,d \mu \in\left[u_{i}-\delta, u_{i}+\delta\right]$ for some $i$. Then
\begin{equation*}
\mathcal{L}(Z) \cap Z_{l}=\bigcup_{m \in \mathbb{N}} \bigcup_{i=1}^{r} B_{m, u_{i}} \tag{4.2}
.\end{equation*}
Denote
$$
D_{N}^{\prime}:=\# \bigcup_{x \in B_{m, i}}\left\{U \in \mathcal{U}^{N}: U\right. \text { satisfies (a) and (b) of Lemma 4.2\}. }
$$
Note that $\mu(Z)=1$. Thus, by Lemma 4.3.3 of [3], there exists $N_{0}>0$ such that
$$
D_{N}^{\prime} \leq \exp \left[N\left(h_{\mu}(\Phi)+2 \delta\right)\right]
$$
holds for all $N \geq N_{0}$.

For all $N \geq N_{0} \in \mathbb{N}$, the collection of all $U \in \mathcal{U}^{N}$ such that Lemma 4.2 holds for some $x \in B_{m, i}$ covers $B_{m, i}$. Denote the set of corresponding ( $x, t$ ) which are from such $U \in \mathcal{U}^{N}$ as
$\Gamma$. Hence, we can estimate
$$
\begin{aligned}
& M\left(B_{m, i}, S, \alpha, \epsilon\right) \\
= & \lim _{T \rightarrow+\infty} \inf _{\Gamma} \sum_{(x, t) \in \Gamma} \exp \left(S_{t}(x, t, \epsilon)-\alpha t\right) \\
\leq & \varlimsup_{N_{0} \rightarrow+\infty} \sum_{N=N_{0}}^{+\infty} D_{N}^{\prime} \exp \left[-\alpha N+N\left(\int_{Z} b d \mu+3 \delta\right)+\gamma_{N}(S, \epsilon)\right] \\
\leq & \varlimsup_{N_{0} \rightarrow+\infty} \sum_{N=N_{0}}^{+\infty} \exp \left[N\left(h_{\mu}(\Phi)+\int_{Z} b d \mu+5 \delta-\alpha+\varlimsup_{t \rightarrow+\infty} \frac{\gamma_{t}(S, \epsilon)}{t}\right)\right] \\
\leq & \varlimsup_{N_{0} \rightarrow+\infty} \sum_{N=N_{0}}^{+\infty} \exp \left[N\left(\sup \left\{h_{\mu}(\Phi)+\int_{Z} b d \mu: \mu \in \mathcal{M}_{Z}\right\}+5 \delta-\alpha+\varlimsup\limits_{t \rightarrow+\infty} \frac{\gamma_{t}(S, \epsilon)}{t}\right)\right]
.\end{aligned}
$$

Let $P:=\sup \left\{h_{\mu}(\Phi)+\int_{Z} b d \mu: \mu \in \mathcal{M}_{Z}\right\}$ and $K:=\exp \left(-\alpha+P+5 \delta+\varlimsup\limits_{t \rightarrow+\infty} \frac{\gamma_{t}(S, \epsilon)}{t}\right)$. Then
$$
M\left(B_{m, i}, S, \alpha, \epsilon\right) \leq \varlimsup_{N_{0} \rightarrow+\infty} \sum_{N=N_{0}}^{+\infty} K^{N}
.$$
Consider $\alpha>P+5 \delta+\varlimsup\limits_{t\rightarrow+\infty}\frac{\gamma_{t}(S, \epsilon)}{t})$. In this case, $K<1$, and we have
$$
M\left(B_{m, i}, S, \alpha, \epsilon\right) \leq \varlimsup_{N_{0} \rightarrow+\infty} \sum_{N=N_{0}}^{+\infty} K^{N}=0.
$$
Hence, $P_{B_{m, i}}(S, \epsilon)<\alpha$. Let $\alpha \rightarrow P+5 \delta+\varlimsup\limits_{t \rightarrow+\infty} \frac{\gamma_{t}(S, \epsilon)}{t}$. Then, we obtain
$$
P_{B_{m, i}}(S, \epsilon) \leq P+5 \delta+\varlimsup\limits _{t \rightarrow+\infty} \frac{\gamma_{t}(S, \epsilon)}{t}
.$$
Combining Lemma 2.3 and equation (4.2), we obtain
$$
\begin{aligned}
P+5 \delta+\lim _{\epsilon \rightarrow 0} \varlimsup_{t \rightarrow+\infty} \frac{\gamma_{t}(S, \epsilon)}{t} & \geq \lim _{\epsilon \rightarrow 0} \sup _{m, i} P_{B_{m, i}}(S, \epsilon) \\
& =\varlimsup_{_{\epsilon \rightarrow 0}} P_{\mathcal{L}(Z) \cap Z_{l}}(S, \epsilon) \\
& =P_{\mathcal{L}(Z) \cap Z_{l}}(S).
\end{aligned}
$$
Since $\varlimsup\limits_{t \rightarrow+\infty} \frac{\gamma_{t}(S, \epsilon)}{t}=0$ (as $S$ satisfies tempered variation), letting $\delta \rightarrow 0$ yields
$$
P_{\mathcal{L}(Z)}(S) \leq P=\sup \left\{h_{\mu}(\Phi)+\int_{Z} b\, d \mu: \mu \in \mathcal{M}_{Z}\right\}
.$$

Finally, we prove that $P_{\mathcal{L}(Z)}(S) \geq \sup \left\{h_{\mu}(\Phi)+\int_{Z} b d \mu: \mu \in \mathcal{M}_{Z}\right\}$.
\begin{lemma}
   ([17]) For each $\epsilon>0$, there exists $0<\delta<\epsilon$, a measurable partition $\xi= \left\{C_{1}, C_{2}, \ldots, C_{m}\right\}$ of $X$, and an open cover $\mathcal{U}=\left\{U_{1}, U_{2}, \ldots, U_{k}\right\}$ of $X$ for some $k \geq m$, such that the following conditions hold:
\begin{enumerate}
  \item $\operatorname{diam}\left(C_{j}\right) \leq \epsilon, \overline{U_{i}} \subseteq C_{i}$, and $\mu\left(C_{i} \backslash U_{i}\right)<\delta$ for $i=1,2, \ldots, m$.
  \item The set $E=\bigcup_{i=m+1}^{k} U_{i}$ has measure $\mu(E)<\delta$.
  \item $2 \delta \log m \leq \epsilon$.
\end{enumerate} 
\end{lemma}
\begin{lemma}
   $P_{Z}(S) \geq h_{\mu}(\Phi)+\int_{Z} b d \mu$ for each $\mu \in \mathcal{M}_{Z}$. 
\end{lemma}
\begin{proof}
  Fix $x \in Z$ and $n \geq 1$. Let $t_{n}(x)$ be the number of $i \in[0, n) \cap \mathbb{N}$ such that $\varphi_{1}^{i}(x) \in E$. By Birkhoff's ergodic theorem, we have
$$
\lim _{n \rightarrow \infty} \frac{t_{n}(x)}{n}=\lim _{n \rightarrow \infty} \frac{1}{n} \sum_{j=0}^{n-1} \chi_{E}\left(\varphi_{1}^{j}(x)\right)=\int_{X} \chi_E\, d \mu=\mu(E)
$$
for $\mu$-almost every $x \in Z$. Then, by Egorov's theorem, there exists $N_{1}>0$ and a subset $Z_{1} \subseteq Z$ such that $\mu\left(Z_{1}\right) \geq 1-\delta$ and
\begin{equation*}
t_{n}(x) \leq 2 \delta n \tag{4.3}
\end{equation*}
for any $x \in Z_{1}$ and $n \geq N_{1}$.\\[0pt]
Lemma 11 of [4] tells us that
$$
\lim _{t \rightarrow \infty} \frac{s_{t}(x)}{t}=\lim _{t \rightarrow \infty} \frac{1}{t} \int_{0}^{t}\left(b \circ \varphi_{u}\right)\, d u=\int_{Z} b \,d \mu
$$
for each $\mu \in \mathcal{M}_{Z}$, for $\mu$-almost every $x \in Z$. Then, by Egorov's theorem, there exists $N_{2}>0$ and a subset $Z_{2} \subseteq Z$ such that $\mu\left(Z_{2}\right) \geq 1-\delta$ and
\begin{equation*}
\left|\frac{s_{n}(x)}{n}-\int_{Z} b d \mu\right|<\delta \tag{4.4}
\end{equation*}
for every $x \in Z_{2}$ and $n \geq N_{2}$.

Let $\xi_{n}=\bigcap_{j=0}^{n} \varphi_{1}^{-j}(\xi)$. Applying the Shannon-McMillan-Breiman theorem and Egorov's theorem, there exists $N_{3} \in \mathbb{N}$ and a measurable set $Z_{3} \subseteq Z$ with $\mu\left(Z_{3}\right) \geq 1-\delta$ such that
\begin{equation*}
\mu\left(\xi_{n}(x)\right) \geq \exp \left[\left(-h_{\mu}\left(\varphi_{1}, \xi\right)+\delta\right) n\right] \tag{4.5}
\end{equation*}  
for every $x \in Z_{3}$ and $n>N_{3}$. Denote $N=\max \left\{N_{1}, N_{2}, N_{3}\right\}$, and let $A=Z_{1} \cap Z_{2} \cap Z_{3}$. We have $\mu(A) \geq 1-3 \delta$, and the inequalities (4.3), (4.4), and (4.5) hold for all $x \in A$ and $n \geq N$.

We can find a large enough $l$ satisfying $\mu\left(Z_{l}\right)>1-\delta$ such that the sets $Z_{l}$ are nested and exhaust $Z$. Thus, $\mu\left(Z_{l}\cap A\right) >1-4 \delta$. Let $\beta$ be the Lebesgue number of the cover $\mathcal{U}$ of $Z_{l}$, and choose $\epsilon_{0}>0$ such that $2 \epsilon_{0}<\beta$. Fix $\alpha \in \mathbb{R}$, and choose $N^{\prime} \geq N$ such that there exists a cover $\Gamma \subseteq X \times[n,+\infty)$ of $Z_{l}$ satisfying
\begin{equation*}
\left|\sum_{(x, t) \in \gamma} \exp \left(S\left(x, t, \epsilon_{0}\right)-\alpha t\right)-M\left(Z_{l}, S, \alpha, \epsilon_{0}\right)\right|<\delta \tag{4.6}
\end{equation*}
holds for each $n \geq N^{\prime}$. Given $q \in \mathbb{N}$, let
$$
\Gamma_{q}=\left\{(x, q) \in \Gamma: B_{q}\left(x, \epsilon_{0}\right) \cap A \neq \emptyset\right\}
$$
and
$$
B_{q}=\bigcup_{(x, t) \in \Gamma_{q}} B_{t}\left(x, \epsilon_{0}\right) .
$$
It follows from Lemma 2 in [17] that
$$
\# \Gamma_{q} \geq \mu\left(B_{q} \cap A\right) \exp \left[h_{\mu}\left(\varphi_{1}, \xi\right) q-(1+2 \log m) q \delta\right]
$$
for all $q \in \mathbb{N}$. Additionally, by equation (4.4), we have
$$
s_{q}(x) \geq q\left(\int_{Z} b\, d \mu-\delta\right)
,$$
and therefore,
$$
\sup _{x \in B_{q}\left(x, \epsilon_{0}\right)} s_{q}(x) \geq q\left(\int_{Z} b\, d \mu-\delta\right)-\gamma_{q}\left(S, \epsilon_{0}\right)
$$
with tempered variation. Thus, for all $q \geq N^{\prime}$ and $(x, t) \in \Gamma_{q}$, we have
$$
\begin{aligned}
& \sum_{(x, t) \in \Gamma} \exp \left(S\left(x, t, \epsilon_{0}\right)-\alpha t\right) \geq \sum_{q=N^{\prime}}^{+\infty} \sum_{(x, t) \in \Gamma_{q}} \exp \left(\sup _{B_{q}\left(x, \epsilon_{0}\right)} s_{q}-\alpha q\right) \\
\geq & \sum_{q=N^{\prime}}^{+\infty} \# \Gamma_{q} \exp \left[\left(-\alpha+\int_{Z} b \,d \mu-\delta\right) q-\gamma_{q}\left(S, \epsilon_{0}\right)\right] \\
\geq & \sum_{q=N^{\prime}}^{+\infty} \mu\left(B_{q} \cap A\right) \exp \left[q\left(h_{\mu}\left(\varphi_{1}, \xi\right)+\int_{Z} b\, d \mu-\frac{\gamma_{q}\left(S, \epsilon_{0}\right)}{q}-\alpha\right)-2(1+\log m) q \delta\right] \\
\geq & \sum_{q=N^{\prime}}^{+\infty} \mu\left(B_{q} \cap A\right) \exp \left[q\left(-\alpha+h_{\mu}\left(\varphi_{1}, \xi\right)+\int_{Z} b \,d \mu-\overline{\lim }_{t \rightarrow+\infty} \frac{\gamma_{t}(S, \epsilon)}{t}-2 \delta(1+\log m)\right)\right] .
\end{aligned}
$$

Let $\alpha<h_{\mu}\left(\varphi_{1}, \xi\right)+\int_{Z} b \,d \mu-\varlimsup\limits_{t \rightarrow+\infty} \frac{\gamma_{t}\left(S, \epsilon_{0}\right)}{t}-2(1+\log m) \delta$. Since $\delta$ can be sufficiently small, we have
$$
\sum_{(x, t) \in \Gamma} \exp \left(S\left(x, t, \epsilon_{0}\right)-\alpha t\right) \geq \sum_{q=N^{\prime}}^{+\infty} \mu\left(B_{q} \cap A\right) \geq \mu(A) \geq 1-3 \delta
.$$
Thus, one has by (4.6)
$$
M\left(Z_{l}, S, \alpha, \epsilon_{0}\right) \geq 1-4 \delta>0
.$$
Therefore, $P_{Z_{l}}\left(S, \epsilon_{0}\right) \geq \alpha$. Let $\alpha \rightarrow h_{\mu}\left(\varphi_{1}, \xi\right)+\int_{Z} b\, d \mu-\varlimsup\limits_{t \rightarrow+\infty} \frac{\gamma_{t}\left(S, \epsilon_{0}\right)}{t}-2(1+\log m) \delta$ and $\delta \rightarrow 0$. For all $\epsilon<\epsilon_{0}$, we obtain
$$
P_{Z_{l}}(S, \epsilon) \geq h_{\mu}\left(\varphi_{1}, \xi\right)+\int_{Z} b\, d \mu-\varlimsup_{\lim _{t \rightarrow+\infty}} \frac{\gamma_{t}\left(S, \epsilon_{0}\right)}{t}
.$$

For the partition $\xi_{q}$ and open cover $\mathcal{U}_{q}$ as mentioned above, let $\epsilon=\frac{1}{q}$. Take $\epsilon_{q}>0$ for each $q>0$ such that $2 \epsilon_{q}<\frac{1}{q}$ is a Lebesgue number of $\mathcal{U}_{q}$. It follows that
$$
\lim _{q \rightarrow+\infty} h_{\mu}\left(\varphi_{1}, \xi_{q}\right)=h_{\mu}\left(\varphi_{1}\right)
$$
since $\operatorname{diam}\left(\xi_{q}\right) \rightarrow 0$ as $q \rightarrow+\infty$. Therefore, we obtain
$$
\begin{aligned}
P_{Z_{l}}(S) & =\lim _{q \rightarrow+\infty} P_{Z_{l}}\left(S, \epsilon_{q}\right) \\
& \geq \lim _{q \rightarrow+\infty}\left(h_{\mu}\left(\varphi_{1}, \xi_{q}\right)+\int_{Z} b\,d \mu-\lim _{q \rightarrow+\infty} \overline{\lim }_{t \rightarrow+\infty} \frac{\gamma_{t}\left(S, \epsilon_{q}\right)}{t}\right) \\
& =h_{\mu}\left(\varphi_{1}\right)+\int_{Z} b \,d \mu
.\end{aligned}
$$
Thus, $P_{Z_{l}}(S) \geq h_{\mu}\left(\varphi_{1}\right)+\int_{Z} b d \mu$. The proof of Lemma 4.4 is completed.
\end{proof}
Since $\mu \in \mathcal{M}_{Z}$ is ergodic, $Z_{\mu}$ is a non-empty $\Phi$-invariant subset of $\mathcal{L}(Z)$ with $\mu\left(Z_{\mu}\right)=1$. Hence,
$$
P_{\mathcal{L}(Z)}(S) \geq h_{\mu}(\Phi)+\int_{Z_{\mu}} b \,d \mu=h_{\mu}(\Phi)+\int_{Z} b\, d \mu.
$$
This completes the proof of Theorem 1.2.
\end{proof}

\end{document}